\documentclass[12pt]{article}

\usepackage{latexsym, amssymb, amsmath, amscd, amsfonts ,texdraw, mathrsfs, amsthm}

\usepackage{fullpage}

\newtheorem{thm}{Theorem}
\newtheorem{lemma}[thm]{Lemma}
\newtheorem{cor}[thm]{Corollary}

\theoremstyle{definition}

\newtheorem{ex}{Example}
\theoremstyle{remark}
\newtheorem{remark}{Remark}

\def\qed{\nopagebreak\hfill{\rule{4pt}{7pt}}}
\def\proof{{\it{Proof.} \hskip 2pt}}

\def\gp{\mathcal{G}_P}

\def\ri{^{}_i}
\def\mR{\mathcal{R}}

\DeclareMathOperator{\pmaj}{pmaj}
\DeclareMathOperator{\Cr}{cr}
\DeclareMathOperator{\inv}{inv}
\DeclareMathOperator{\maj}{maj}
\DeclareMathOperator{\des}{des}

\begin{document}

\title{A Major Index for Matchings and Set Partitions}

\author{
William Y.\,C. Chen$^{1,5}$, Ira M. Gessel$^{2}$,
Catherine H. Yan$^{3,6}$ and \\
Arthur L.\,B. Yang$^{4,5}$
 \vspace{.2cm} \\
$^{1,3,4}$Center for Combinatorics, LPMC\\
Nankai University, Tianjin 300071, P. R. China
\vspace{.2cm} \\
$^{2}$Department of Mathematics \\
Brandeis University, Waltham, MA 02454-9110  \vspace{.2cm} \\
$^3$Department of Mathematics\\
Texas A\&M University, College Station, TX 77843 \vspace{.2cm} \\
{\small E-mail: $^1$chen@nankai.edu.cn, $^2$gessel@brandeis.edu,
$^3$cyan@math.tamu.edu, $^4$yang@nankai.edu.cn}}

\date{}
\maketitle

\begin{abstract}
 We introduce a
statistic $\pmaj$ on  partitions of $[n]=\{1,2,\dots, n\}$, and show that it is
equidistributed with the number of 2-crossings over
partitions of $[n]$  with given sets of minimal block elements and
maximal block elements.
This  generalizes the classical result of equidistribution for the
permutation statistics  inversion number and major index.
\end{abstract}



\footnotetext[5]{The first and the fourth authors were
supported by the 973 Project on
Mathematical Mechanization, the Ministry of Education, the Ministry
of Science and Technology, and the National Science Foundation of
China.}

\footnotetext[6]{The third author was supported in part by NSF
grant \#DMS-0245526.}

\section{Introduction}
One of the classical results on permutations is the equidistribution
of the statistics inversion number and major index. For a permutation $\pi=(a_1a_2\cdots
a_n)$, a pair $(a_i, a_j)$ is called an
\emph{inversion} if $i < j$ and $a_i > a_j$.  The statistic  $\inv(\pi)$ is
defined as the
number of inversions of $\pi$. The \emph{descent set}
$D(\pi)$ is defined as $\{\,i: a_i> a_{i+1}\,\}$ and its cardinality is
denoted by $\des(\pi)$.
The sum of the elements of $D(\pi)$ is called the \emph{major index}
of $\pi$ (also called the greater index) and
denoted by $\maj(\pi)$.
Similarly one can define the notions of inversion, descent set, and
major index for any word $w=w_1w_2\cdots w_n$ of not necessarily
distinct integers.
It is a  result of MacMahon \cite{MacMahon}
that $\inv$ and $\maj$ are
equidistributed on the rearrangement class of any word.
A statistic equidistributed with $\inv$ is called \emph{Mahonian}.

There are many research articles devoted to Mahonian statistics and
their generalizations. For  example, see \cite{clarke,Han} for
Mahonian statistics for words, \cite{Si-St94,Si-St96} for Mahonian
statistics and Laguerre polynomials,  \cite{Sagan} for a major index
statistic for  set partitions, and very recently \cite{HS06} for
$\inv$ and $\maj$ for standard Young tableaux.

Given a partition of $[n]$, there is a natural generalization of
inversions, namely, 2-crossings, which can be viewed easily on a
graphical representation of the partition. In this paper we
introduce a new statistic, called the \emph{p-major index} and
denoted $\pmaj(P)$, on the set of partitions of $[n]$. We prove that
for any $S, T \subseteq [n]$ with $|S|=|T|$, $\pmaj$ and $\Cr_2$,
the number of 2-crossings,  are equally distributed on the set
$P_n(S, T)$. Here $P_n(S, T)$ is the set of  partitions of $[n]$ for
which $S$ is the set of minimal block elements, and $T$ is the set
of maximal block elements. Restricted to permutations, the pair
$(\Cr_2, \pmaj)$ coincides with $(\inv, \maj)$.  Hence our result
gives another generalization of MacMahon's equidistribution theorem.

In the next section we introduce  notation and state the main
results. An algebraic proof and some examples are given in Section
3. In Section 4 we present a bijective proof which generalizes
Foata's second fundamental transformation \cite{Foata,Foata97}.

\section{Definitions and the main results}
A {\em partition} of $[n]=\{1,2,\dots, n\}$ is a collection $P$ of disjoint
nonempty subsets of $[n]$, whose union is $[n]$. Each subset in $P$
is called a
{\em block} of $P$. A {\em (perfect) matching}
of $[n]$ is a partition of $[n]$ in which each block contains exactly
two elements. We denote by $\Pi_n$ the set of all partitions of $[n]$.
Following \cite{CDDSY},
we represent each  partition $P\in \Pi_n$ by a graph $\gp$ on the
vertex set $[n]$ whose edge set
consists of arcs connecting the elements of each block in numerical
order.  Such a graph is called the \emph{standard representation}
of the partition $P$.  For example, the standard representation of
1457-26-3 has the arc set $\{(1,4), (4, 5), (5, 7), (2,6)\}$.
We always write
an arc $e$ as a pair $(i,j)$ with $i < j$, and say that $i$ is the
\emph{left-hand endpoint} of $e$ and $j$ is the \emph{right-hand
endpoint} of $e$.  \\

\begin{figure}[ht]
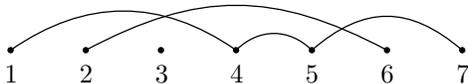

\centertexdraw{
\drawdim mm

\move(10 7) \fcir f:0 r:0.4
\move(20 7) \fcir f:0 r:0.4
\move(30 7) \fcir f:0 r:0.4
\move(40 7) \fcir f:0 r:0.4
\move(50 7) \fcir f:0 r:0.4
\move(60 7) \fcir f:0 r:0.4
\move(70 7) \fcir f:0 r:0.4

\move(11 4) \textref h:R v:C \htext{\footnotesize{$1$}}
\move(21 4) \textref h:R v:C \htext{\footnotesize{$2$}}
\move(31 4) \textref h:R v:C \htext{\footnotesize{$3$}}
\move(41 4) \textref h:R v:C \htext{\footnotesize{$4$}}
\move(51 4) \textref h:R v:C \htext{\footnotesize{$5$}}
\move(61 4) \textref h:R v:C \htext{\footnotesize{$6$}}
\move(71 4) \textref h:R v:C \htext{\footnotesize{$7$}}

\linewd 0.2
\move(10 7) \clvec (20 14)(30 14)(40 7)
\move(40 7) \clvec (43 10)(47 10)(50 7)
\move(50 7) \clvec (57 13)(63 13)(70 7)
\move(20 7) \clvec (35 15)(45 15)(60 7)
}
\caption{The standard representation of partition $P=1457-26-3$.}
\label{partition}
\end{figure}

A  partition $P \in \Pi_n$ is a matching if and
only if in $\gp$ each vertex is the endpoint of exactly one arc. In
other words, each vertex  is either a left-hand endpoint or a
right-hand endpoint. In particular, a permutation $\pi$ of $[m]$ can
be represented  as a matching $M_\pi$ of $[2m]$ with arcs connecting
$m+1-\pi(i)$ and $i+m$ for $ 1\leq i \leq m$. See Figure 2 for an
example.

\vskip 4mm

\begin{figure}[ht]
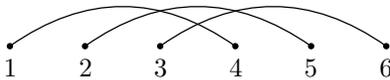

\centertexdraw{
\drawdim mm

\linewd 0.2
\move(10 10) \clvec (20 17)(30 17)(40 10)
\move(20 10) \clvec (30 17)(40 17)(50 10)
\move(30 10) \clvec (40 17)(50 17)(60 10)

\move(10 10) \fcir f:0 r:0.4
\move(20 10) \fcir f:0 r:0.4
\move(30 10) \fcir f:0 r:0.4
\move(40 10) \fcir f:0 r:0.4
\move(50 10) \fcir f:0 r:0.4
\move(60 10) \fcir f:0 r:0.4

\move(11 7) \textref h:R v:C \htext{\footnotesize{$1$}}
\move(21 7) \textref h:R v:C \htext{\footnotesize{$2$}}
\move(31 7) \textref h:R v:C \htext{\footnotesize{$3$}}
\move(41 7) \textref h:R v:C \htext{\footnotesize{$4$}}
\move(51 7) \textref h:R v:C \htext{\footnotesize{$5$}}
\move(61 7) \textref h:R v:C \htext{\footnotesize{$6$}}
}
\caption{The permutation $\pi=321$ and the matching $M_\pi$.}
\label{f1}
\end{figure}

Two arcs $(i_1, j_1), (i_2, j_2)$ of $\gp$ form  a 2-crossing if
$i_1 < i_2 < j_1 < j_2$. Let $\Cr_2(P)$ denote the number of
2-crossings of $P$.  A 2-crossing is a natural generalization of an inversion
of a permutation. It is easily seen that
under the correspondence $\pi \mapsto M_\pi$,
 $\Cr_2(M_\pi)=\inv(\pi)$.

Given $P \in \Pi_n$, define
\begin{eqnarray*}
\min(P)=\{ \text{minimal block elements of $P$}\}, \\
\max(P)=\{ \text{maximal block elements of $P$}\}.
\end{eqnarray*}
For example, for $P=\mbox{1457-26-3}$, $\min(P) =\{1,2,3\}$ and
$\max(P)=\{3,6,7\}$.

Fix $S, T \subseteq [n]$ with $|S|=|T|$.  Let $P_n(S, T)$ be the
set $\{\, P \in \Pi_n: \min(P)=S, \max(P)=T\,\}$.
For any set $X \subseteq [n]$, let  $X\ri =X \cap \{i+1, \dots, n\}$.
\begin{thm} \label{2crossing}
Fix $S, T \subseteq [n]$ with $|S|=|T|$.
Then
\begin{eqnarray}
\sum_{P \in P_n(S, T)} y^{\Cr_2(P)} = \prod_{i \notin T} (1+y+\cdots
+y^{h(i)-1}),
\end{eqnarray}
where  $h(i)=|T\ri|-|S\ri|$.
\end{thm}

For a permutation $\pi=(a_1a_2\cdots a_n)$, the major index
$\maj(\pi)$ can be computed as $\sum_{i=1}^n \des(a_i\cdots a_n)$.
This motivates the following definition of the p-major index for set
partitions. Given $P \in \Pi_n$, we start with the standard
representation $\gp$. First label the arcs of $P$ by $1, 2, \dots,
k$ from right to left in order of their left-hand endpoints. That is,
if the arcs are $(i_1, j_1), (i_2, j_2), \dots, (i_k, j_k)$ with
$i_1 > i_2 > \cdots > i_k$, then $(i_r, j_r)$ has label $r$, for $1
\leq r \leq k$. Next we associate a sequence $\sigma(r)$ to each
right-hand endpoint $r$. Assume that the right-hand endpoints are $r_1
> r_2 > \cdots > r_k$. (The set $\{r_1, \dots, r_k\}$ is exactly
$[n]\setminus S$.) The sequence $\sigma(r_i)$ is defined backward
recursively: let $\sigma(r_1)=a$ if  $r_1$ is the right-hand endpoint
of the arc with label $a$. In general, after defining $\sigma(r_i)$,
assume that the left-hand endpoints of the arcs labeled $a_1, \dots,
a_t$ lie between $r_{i+1}$ and $r_i$. Then $\sigma(r_{i+1})$
is obtained from $\sigma(r_i)$ by deleting entries $a_1, \dots, a_t$
and adding $b$ at the very beginning, where $b$ is the label for the
arc whose right-hand endpoint is $r_{i+1}$. Finally, define the
statistic $\pmaj(P)$ by
\[
\pmaj(P):= \sum_{r_i} \des(\sigma(r_i)).
\]

\begin{ex}
\begin{figure}[ht]
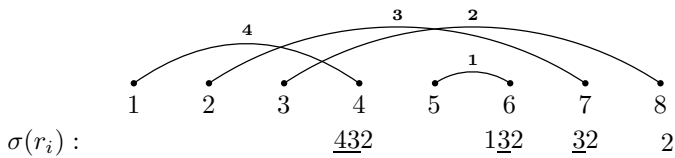

\centertexdraw{
\drawdim mm

\linewd 0.2
\move(10 10) \clvec (20 17)(30 17)(40 10)
\move(20 10) \clvec (35 20)(55 20)(70 10)
\move(30 10) \clvec (45 20)(65 20)(80 10)
\move(50 10) \clvec (53 12)(57 12)(60 10)

\move(56 13) \textref h:R v:C \htext{\tiny{\bf 1}}
\move(56 19) \textref h:R v:C \htext{\tiny{\bf 2}}
\move(46 19) \textref h:R v:C \htext{\tiny{\bf 3}}
\move(26 17) \textref h:R v:C \htext{\tiny{\bf 4}}

\move(10 10) \fcir f:0 r:0.4
\move(20 10) \fcir f:0 r:0.4
\move(30 10) \fcir f:0 r:0.4
\move(40 10) \fcir f:0 r:0.4
\move(50 10) \fcir f:0 r:0.4
\move(60 10) \fcir f:0 r:0.4
\move(70 10) \fcir f:0 r:0.4
\move(80 10) \fcir f:0 r:0.4

\move(11 7) \textref h:R v:C \htext{\footnotesize{$1$}}
\move(21 7) \textref h:R v:C \htext{\footnotesize{$2$}}
\move(31 7) \textref h:R v:C \htext{\footnotesize{$3$}}
\move(41 7) \textref h:R v:C \htext{\footnotesize{$4$}}
\move(51 7) \textref h:R v:C \htext{\footnotesize{$5$}}
\move(61 7) \textref h:R v:C \htext{\footnotesize{$6$}}
\move(71 7) \textref h:R v:C \htext{\footnotesize{$7$}}
\move(81 7) \textref h:R v:C \htext{\footnotesize{$8$}}

\move(3 2) \textref h:R v:C \htext{\footnotesize{$\sigma(r_i):$}}
\move(42 2) \textref h:R v:C \htext{\footnotesize{$\underline{4}
    \underline{3}2$}}
\move(62 2) \textref h:R v:C \htext{\footnotesize{$1\underline{3}2$}}
\move(72 2) \textref h:R v:C \htext{\footnotesize{$\underline{3}2$}}
\move(82 2) \textref h:R v:C \htext{\footnotesize{$2$}}

}
\caption{The major index for partition $14-27-38-56$ is 4.}
\label{major_f1}
\end{figure}

Let $P=14-27-38-56$. Then $(r_1,r_2,r_3,r_4)=(8, 7, 6, 4)$.
Figure \ref{major_f1} shows how to compute $\pmaj(P)$.
The sequences $\sigma(r_i)$ are $\sigma(8)=(2)$,  $\sigma(7)=(32)$,
$\sigma(6)=(132)$, and $\sigma(4)=(432)$. For each $\sigma(r_i)$,
the elements in the  descent set are underlined.
The p-major index of $P$ is $\pmaj(P)= 1+1+2=4$.

\end{ex}

\begin{thm} \label{major}
Fix $S, T \subseteq [n]$ with $|S|=|T|$. Then
\begin{eqnarray} \label{for-major}
\sum_{P \in P_n(S, T)} y^{\pmaj(P)} = \prod_{i \notin T} (1+y+\cdots
+y^{h(i)-1}),
\end{eqnarray}
where   $h(i)=|T\ri|-|S\ri|$.
\end{thm}

Combining Theorems \ref{2crossing} and \ref{major}, we have
\begin{cor}\label{cr=maj}
For each $P_n(S, T)$,
$$
\sum_{P \in P_n(S, T)} y^{\pmaj(P)}=\sum_{P \in P_n(S, T)} y^{\Cr_2(P)}.
$$
That is, the two statistics $\Cr_2$ and $\pmaj$ have the same
distribution over each set $P_n(S, T)$.
\end{cor}

When $n=2m$, $S=[m]$ and $T=[2m]\setminus [m]$,
the map $\pi \mapsto M_\pi$ gives a one-to-one correspondence
between $P_n(S, T)$ and the set of permutations of $[m]$.
It is easy to see that  $\pmaj(M_\pi)=\maj(\pi)$.
Hence the equidistribution of $\inv$ and $\maj$ for permutations is a
special case of Corollary \ref{cr=maj}.

Another consequence of Theorems \ref{2crossing} and \ref{major} is
the symmetry of the number of partitions of $[n]$ with a given
number of 2-crossings or a given p-major index. Let $A(n,i; S,
T)$ be the set of partitions in $P_n(S, T)$ such that $\Cr_2(P)=i$,
whose cardinality is $a(n,i;S, T)$.
\begin{cor}
Fix $n$ and let $K=\sum_{i \notin T} (h(i)-1)$.
Then the sequence $\{ a(n,i; S, T)\}_{i=0}^{K}$ is symmetric. That is,
$$
a(n,i; S, T) = a(n, K-i; S, T).
$$
The same result holds if we replace $\Cr_2(P)$ by $\pmaj(P)$ in defining
$A(n,i;S,T)$ and $a(n,i;S,T)$.
\end{cor}

\section{Proofs for the  main results}

In this section we give the proofs for Theorems \ref{2crossing} and
\ref{major}. Given a partition $P \in \Pi_n$, a vertex $i \in [n]$
in the standard representation $\gp$ is  one of the following types:
\begin{enumerate}
\item a left-hand endpoint if $i \in \min(P)\setminus \max(P)$,
\item a right-hand endpoint if $i \in \max(P) \setminus  \min(P)$,
\item an isolated point if $i \in \min(P) \cap \max(P)$,
\item a left-hand endpoint and a right-hand endpoint if $i \notin
  \min(P)\cup \max(P)$.
\end{enumerate}
In particular, $[n]\setminus \max(P)$ is the set of points which are
the left-hand endpoints of some arcs, and $[n]\setminus \min(P)$ is
the set of right-hand endpoints.   Fixing $\min(P)=S$ and
$\max(P)=T$ is equivalent to fixing the type of each vertex in
$[n]$. Since the standard representation uniquely determines the
partition, we can identify a  partition $P \in \Pi_n$ with the set
of arcs of $\gp$. Hence the set
 $P_n(S,T)$ is in one-to-one correspondence with the set of matchings
between $[n]\setminus T$ and $[n]\setminus S$ such that
$i < j$ whenever  $i \in [n] \setminus T$, $j \in
[n]\setminus S$ and $i$ is matched to $j$.
In the following such a matching is referred as a \emph{good}
matching. Denote by $M_n(S, T)$ the set of good matchings
from  $[n]\setminus T$ to $[n]\setminus S$.

\medskip

\noindent \textbf{Proof of Theorem \ref{2crossing}}.
Assume $[n]\setminus T=\{i_1, i_2, \dots, i_k\}$ with $i_1 < i_2 <
\dots < i_k$.
Let $S(H)$ be the set of sequences $\{(a_1, a_2, \dots, a_k)\}$ where
$1 \leq a_r \leq h(i_r)$ for each $1 \leq r \leq k$. We give a
bijection between the sets $M_n(S,T)$ and $S(H)$.
The construction is essentially due to M.~de Sainte-Catherine \cite{SC83}.

Given a sequence $\alpha=(a_1, a_2, \dots, a_k)$ in $S(H)$, we
construct a matching from $[n]\setminus T=\{i_1, i_2, \dots, i_k\}$
to $[n]\setminus S$ as follows. First, there are exactly $h(i_k)$
many elements in $[n] \setminus S$ which are greater than $i_k$.
List them in increasing order as $1, 2,  \dots, h(i_k)$. Match $i_k$
to the $a_k$th element, and mark this element as  \emph{dead}.

In general, after matching elements $i_{j+1}, \dots, i_k$ to some
elements in $[n]\setminus S$, we process the element $i_j$. At this
stage there are exactly $h(i_j)$ many elements in $[n]\setminus S$
which are greater than $i_j$ and not dead. List them in increasing
order by $1, 2, \dots, h(i_j)$. Match $i_j$ to the $a_j$th of them,
and mark it as dead. Continuing the process  until $j=1$, we get a
good matching $M(\alpha) \in M_n(S,T)$. The map $\alpha \mapsto
M(\alpha)$
 gives the desired
bijection between  $S(H)$ and $M_n(S, T)$.

Let $P(\alpha)$ be the partition of $[n]$ for which the arc set of $\gp$ is
$M(\alpha)$.
By the above construction, the number of 2-crossings formed by arcs
$(i_j, b)$ and $(a,c)$ with $a < i_j < c < b$ is exactly $a_j-1$.
Hence $\Cr_2(P(\alpha))=\sum_{j=1}^k (a_j-1)$ and
\begin{eqnarray}
\sum_{P \in P_n(S, T)} y^{\Cr_2(P)}=\sum_{(a_1, \dots, a_k) \in S(H)}
y^{\sum_{j=1}^k(a_j-1)} = \prod_{i \notin T} (1+y+\cdots +y^{h(i)-1}).
\end{eqnarray}
\qed

\begin{ex} \label{ex-crossing}
Let $n=6$, $S=\{1,2\}$, and $T=\{5,6\}$.
Then $[n]\setminus T=\{1,2,3,4\}$, $[n]\setminus S=\{3,4,5,6\}$,
and $h(1)=1$, $h(2)=h(3)=h(4)=2$. Figure \ref{crossing-ex} shows the
correspondence between $S(H)$ and $P_n(S, T)$.
For simplicity we omit the vertex labeling.
\begin{figure}[ht]
\centertexdraw{
\drawdim mm
\linewd 0.2
\move(10 100)\textref h:R v:C \htext{{sequence}}
\move(60 100) \textref h:R v:C \htext{{partition}}
\move(102 100) \textref h:R v:C \htext{{$\Cr_2(P)$}}

\move(10 90)\textref h:R v:C \htext{\footnotesize{$(1,1,1,1)$}}
\move(10 80)\textref h:R v:C \htext{\footnotesize{$(1,1,1,2)$}}
\move(10 70)\textref h:R v:C \htext{\footnotesize{$(1,1,2,1)$}}
\move(10 60)\textref h:R v:C \htext{\footnotesize{$(1,1,2,2)$}}
\move(10 50)\textref h:R v:C \htext{\footnotesize{$(1,2,1,1)$}}
\move(10 40)\textref h:R v:C \htext{\footnotesize{$(1,2,1,2)$}}
\move(10 30)\textref h:R v:C \htext{\footnotesize{$(1,2,2,1)$}}
\move(10 20)\textref h:R v:C \htext{\footnotesize{$(1,2,2,2)$}}

\move(100 90)\textref h:R v:C \htext{\footnotesize{$0$}}
\move(100 80)\textref h:R v:C \htext{\footnotesize{$1$}}
\move(100 70)\textref h:R v:C \htext{\footnotesize{$1$}}
\move(100 60)\textref h:R v:C \htext{\footnotesize{$2$}}
\move(100 50)\textref h:R v:C \htext{\footnotesize{$1$}}
\move(100 40)\textref h:R v:C \htext{\footnotesize{$2$}}
\move(100 30)\textref h:R v:C \htext{\footnotesize{$2$}}
\move(100 20)\textref h:R v:C \htext{\footnotesize{$3$}}

\move(30 90) \fcir f:0 r:0.4
\move(40 90) \fcir f:0 r:0.4
\move(50 90) \fcir f:0 r:0.4
\move(60 90) \fcir f:0 r:0.4
\move(70 90) \fcir f:0 r:0.4
\move(80 90) \fcir f:0 r:0.4

\move(30 90) \clvec (50  98 )(60 98)(80 90)
\move(40 90) \clvec (43  95 )(47 95)(50 90)
\move(50 90) \clvec (53  95 )(57 95)(60 90)
\move(60 90) \clvec (63  95 )(67 95)(70 90)

\move(30 80) \fcir f:0 r:0.4
\move(40 80) \fcir f:0 r:0.4
\move(50 80) \fcir f:0 r:0.4
\move(60 80) \fcir f:0 r:0.4
\move(70 80) \fcir f:0 r:0.4
\move(80 80) \fcir f:0 r:0.4

\move(30 80) \clvec (45  87 )(55 87)(70 80)
\move(40 80) \clvec (43  85 )(47 85)(50 80)
\move(50 80) \clvec (53  85 )(57 85)(60 80)
\move(60 80) \clvec (67  86 )(73 86)(80 80)

\move(30 70) \fcir f:0 r:0.4
\move(40 70) \fcir f:0 r:0.4
\move(50 70) \fcir f:0 r:0.4
\move(60 70) \fcir f:0 r:0.4
\move(70 70) \fcir f:0 r:0.4
\move(80 70) \fcir f:0 r:0.4

\move(30 70) \clvec (42  77 )(48 77)(60 70)
\move(40 70) \clvec (43  75 )(47 75)(50 70)
\move(50 70) \clvec (62  77 )(68 77)(80 70)
\move(60 70) \clvec (63  75 )(67 75)(70 70)

\move(30 60) \fcir f:0 r:0.4
\move(40 60) \fcir f:0 r:0.4
\move(50 60) \fcir f:0 r:0.4
\move(60 60) \fcir f:0 r:0.4
\move(70 60) \fcir f:0 r:0.4
\move(80 60) \fcir f:0 r:0.4

\move(30 60) \clvec (42  67 )(48 67)(60 60)
\move(40 60) \clvec (43  65 )(47 65)(50 60)
\move(50 60) \clvec (58  66 )(62 66)(70 60)
\move(60 60) \clvec (68  66 )(72 66)(80 60)

\move(30 50) \fcir f:0 r:0.4
\move(40 50) \fcir f:0 r:0.4
\move(50 50) \fcir f:0 r:0.4
\move(60 50) \fcir f:0 r:0.4
\move(70 50) \fcir f:0 r:0.4
\move(80 50) \fcir f:0 r:0.4

\move(30 50) \clvec (38  56 )(42 56)(50 50)
\move(40 50) \clvec (55  57 )(65 57)(80 50)
\move(50 50) \clvec (53  55 )(57 55)(60 50)
\move(60 50) \clvec (63  55 )(67 55)(70 50)

\move(30 40) \fcir f:0 r:0.4
\move(40 40) \fcir f:0 r:0.4
\move(50 40) \fcir f:0 r:0.4
\move(60 40) \fcir f:0 r:0.4
\move(70 40) \fcir f:0 r:0.4
\move(80 40) \fcir f:0 r:0.4

\move(30 40) \clvec (38  46 )(42 46)(50 40)
\move(40 40) \clvec (52  47 )(57 47)(70 40)
\move(50 40) \clvec (53  45 )(57 45)(60 40)
\move(60 40) \clvec (68  46 )(72 46)(80 40)

\move(30 30) \fcir f:0 r:0.4
\move(40 30) \fcir f:0 r:0.4
\move(50 30) \fcir f:0 r:0.4
\move(60 30) \fcir f:0 r:0.4
\move(70 30) \fcir f:0 r:0.4
\move(80 30) \fcir f:0 r:0.4

\move(30 30) \clvec (38  36 )(42 36)(50 30)
\move(40 30) \clvec (48  36 )(52 36)(60 30)
\move(50 30) \clvec (62  37 )(68 37)(80 30)
\move(60 30) \clvec (63  35 )(67 35)(70 30)

\move(30 20) \fcir f:0 r:0.4
\move(40 20) \fcir f:0 r:0.4
\move(50 20) \fcir f:0 r:0.4
\move(60 20) \fcir f:0 r:0.4
\move(70 20) \fcir f:0 r:0.4
\move(80 20) \fcir f:0 r:0.4

\move(30 20) \clvec (37  25 )(43 25)(50 20)
\move(40 20) \clvec (47  25 )(53 25)(60 20)
\move(50 20) \clvec (57  25 )(63 25)(70 20)
\move(60 20) \clvec (67  25 )(73 25)(80 20)
}
\caption{The $\Cr_2(P)$ for $P \in P_6(\{1,2\}, \{5,6\})$. }
\label{crossing-ex}
\end{figure}
\end{ex}

To prove Theorem \ref{major}, we need a lemma on permutations.
\begin{lemma} \label{lemma}
Let $\sigma=a_1a_2\cdots a_{n-1}$ be a permutation of $\{2,3,\dots,
n\}$. Let $\sigma_0=1a_1\cdots a_{n-1}$ and
$\sigma_i$ be obtained from $\sigma$ by inserting
$1$ right after $a_i$.
Then of the $n$ permutations $\sigma_0, \dots, \sigma_{n-1}$, the major
indices are all distinct and run from $\maj(\sigma)$ to
$\maj(\sigma)+n-1$ in some order.
\end{lemma}
\noindent {\bf Proof.} First note that
$\maj(\sigma_0)=\maj(\sigma)+\des(\sigma)$. Assume that there are $t_i$
descents of $\sigma$ that are greater than $i$. Then
\begin{eqnarray*}
\maj(\sigma_i)=\left\{ \begin{array}{ll}
     \maj(\sigma)+t_i     & \text{ if } a_i > a_{i+1}, \\
     \maj(\sigma)+i+t_i   & \text{ if } a_i <a_{i+1} \text{ or } i=n-1.
     \end{array} \right.
\end{eqnarray*}
It can be checked  that the major indices
of $\sigma_0, \dots, \sigma_{n-1}$
are all distinct and run from $\maj(\sigma)$ to
$\maj(\sigma)+n-1$ in some order.
\qed

Lemma \ref{lemma} is a special case of a result of Stanley \cite[equation (24)]{Stanley},
which implies that the same conclusion holds when inserting any number not occurring in
$\sigma$ into all possible positions in $\sigma$.
The case in which  $n$ is inserted instead of $1$
was used by Gupta \cite{Gupta} to get the generating function
of the  major index over all permutations of $[n]$:
\begin{eqnarray}\label{permu}
\sum_{\pi \in S_n} y^{\maj(\pi)} = \frac{(1-y)(1-y^2)\cdots
  (1-y^n)}{(1-y)^n}.
\end{eqnarray}
This formula  is the special case of Theorem \ref{major}
 with $S=[n]$ and $T=[2n]\setminus [n]$.

\medskip

\noindent {\bf Proof of Theorem \ref{major}.}
Consider the contribution of the arc with label 1 to the
generating function $\sum_{P \in P_n(S,T)} y^{\pmaj(P)}$.
Again we identify the set $P_n(S, T)$ with the set $M_n(S, T)$ of
good matchings from $[n]\setminus T$ to $[n]\setminus S$.

Let $i_k=\max([n]\setminus T)$, which is the left-hand endpoint of
the arc labeled by $1$ in the definition of $\pmaj(P)$, for any $P
\in P_n(S,T)$. Assume $T^{(r)}_{i_k}=\{j_1, j_2, \dots,
j_{h(i_k)}\}$. Let $A=[n]\setminus (T \cup \{i_k\})$, and
$B=[n]\setminus (S \cup \{j_{h(i_k)}\})$. For any good matching $M$
between $A$ and $B$ let $M_t$ ($ 1 \leq t \leq h(i_k)$) be the
matching obtained from $M$ by joining the pair $(i_k, j_t)$, and
replacing each pair $(a, j_r)$, for $r > t$, with  $(a, j_{r+1})$.
Consequently, the arc labeling of $M_t$ can be obtained from that of
$M$ by labeling the arc $(i_k, j_t)$ by $1$, and adding $1$ to the
label  of each arc of $M$. Assume $\sigma(j_1)=b_1b_2\dots
b_{h(i_k)-1}$ for $M$. Then by the definition of $\pmaj$, we have
$$
\pmaj(M_t) = \pmaj(M) + \maj(b'_1\cdots b'_{t-1} 1 b'_t \cdots
b'_{h(i_k)-1})-\maj(b_1\cdots b_{h(i_k)-1}),
$$
where $b'_i=b_i+1$. By Lemma \ref{lemma},  the values of
$\maj(b'_1\cdots b'_{t-1} 1 b'_t \cdots
b'_{h(i_k)-1})-\maj(b_1\cdots b_{h(i_k)-1})$ are all distinct and
run over the set $\{0,1, \dots, h(i_k)-1\}$. Hence
$$
\sum_{P \in P_n(S,T)} y^{\pmaj(P)}=(1+y+\cdots +y^{h(i_k)-1}) \sum_{P
  \in P_{n-1}(A, B)} y^{\pmaj(P)}.
$$
Then equation \eqref{for-major} follows by induction.
\qed

\begin{ex}
The p-major indices for the partitions in Example \ref{ex-crossing}
are given in Figure \ref{major-ex}.  For simplicity we omit the vertex
labeling, but put the sequence $\sigma(r)$ under each right-hand
endpoint $r$.

\begin{figure}[ht]
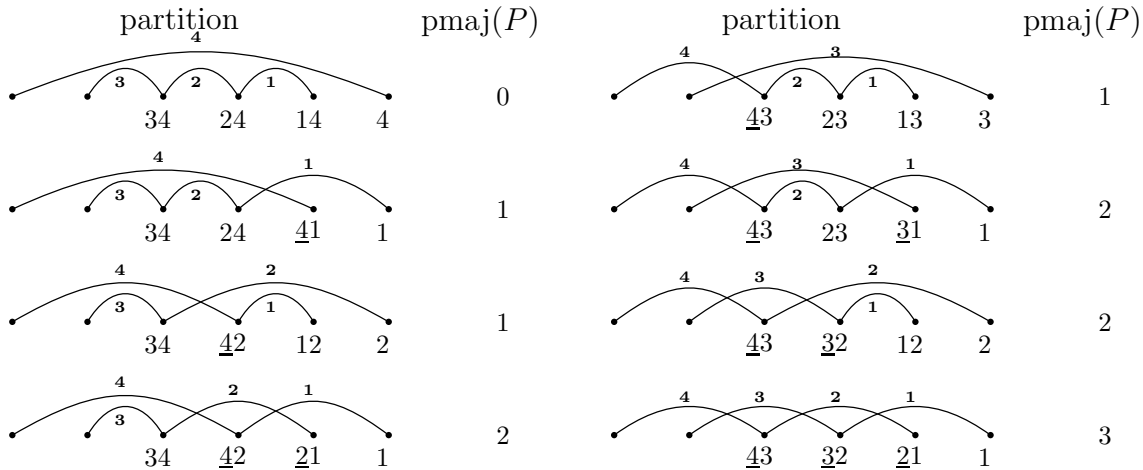

\centertexdraw{
\drawdim mm
\linewd 0.2
\move(60 100) \textref h:R v:C \htext{{partition}}
\move(100 100) \textref h:R v:C \htext{{$\pmaj(P)$}}
\move(140 100) \textref h:R v:C \htext{{partition}}
\move(180 100) \textref h:R v:C \htext{{$\pmaj(P)$}}

\move(96 90)\textref h:R v:C \htext{\footnotesize{$0$}}
\move(96 75)\textref h:R v:C \htext{\footnotesize{$1$}}
\move(96 60)\textref h:R v:C \htext{\footnotesize{$1$}}
\move(96 45)\textref h:R v:C \htext{\footnotesize{$2$}}
\move(176 90)\textref h:R v:C \htext{\footnotesize{$1$}}
\move(176 75)\textref h:R v:C \htext{\footnotesize{$2$}}
\move(176 60)\textref h:R v:C \htext{\footnotesize{$2$}}
\move(176 45)\textref h:R v:C \htext{\footnotesize{$3$}}

\move(30 90) \fcir f:0 r:0.4
\move(40 90) \fcir f:0 r:0.4
\move(50 90) \fcir f:0 r:0.4
\move(60 90) \fcir f:0 r:0.4
\move(70 90) \fcir f:0 r:0.4
\move(80 90) \fcir f:0 r:0.4

\move(30 90) \clvec (50  98 )(60 98)(80 90)
\move(40 90) \clvec (43  95 )(47 95)(50 90)
\move(50 90) \clvec (53  95 )(57 95)(60 90)
\move(60 90) \clvec (63  95 )(67 95)(70 90)

\move(65 92) \textref h:R v:C \htext{\tiny{\bf 1}}
\move(55 92) \textref h:R v:C \htext{\tiny{\bf 2}}
\move(45 92) \textref h:R v:C \htext{\tiny{\bf 3}}
\move(55 98) \textref h:R v:C \htext{\tiny{\bf 4}}

\move(80 87)\textref h:R v:C \htext{\footnotesize{4}}
\move(71 87)\textref h:R v:C \htext{\footnotesize{14}}
\move(61 87)\textref h:R v:C \htext{\footnotesize{24}}
\move(51 87)\textref h:R v:C \htext{\footnotesize{34}}

\move(30 75) \fcir f:0 r:0.4
\move(40 75) \fcir f:0 r:0.4
\move(50 75) \fcir f:0 r:0.4
\move(60 75) \fcir f:0 r:0.4
\move(70 75) \fcir f:0 r:0.4
\move(80 75) \fcir f:0 r:0.4

\move(30 75) \clvec (45  82 )(55 82)(70 75)
\move(40 75) \clvec (43  80 )(47 80)(50 75)
\move(50 75) \clvec (53  80 )(57 80)(60 75)
\move(60 75) \clvec (67  81 )(73 81)(80 75)

\move(70 81) \textref h:R v:C \htext{\tiny{\bf 1}}
\move(55 77) \textref h:R v:C \htext{\tiny{\bf 2}}
\move(45 77) \textref h:R v:C \htext{\tiny{\bf 3}}
\move(50 82) \textref h:R v:C \htext{\tiny{\bf 4}}

\move(80 72)\textref h:R v:C \htext{\footnotesize{1}}
\move(71 72)\textref h:R v:C \htext{\footnotesize{\underline{4}1}}
\move(61 72)\textref h:R v:C \htext{\footnotesize{24}}
\move(51 72)\textref h:R v:C \htext{\footnotesize{34}}

\move(30 60) \fcir f:0 r:0.4
\move(40 60) \fcir f:0 r:0.4
\move(50 60) \fcir f:0 r:0.4
\move(60 60) \fcir f:0 r:0.4
\move(70 60) \fcir f:0 r:0.4
\move(80 60) \fcir f:0 r:0.4

\move(30 60) \clvec (42  67 )(48 67)(60 60)
\move(40 60) \clvec (43  65 )(47 65)(50 60)
\move(50 60) \clvec (62  67 )(68 67)(80 60)
\move(60 60) \clvec (63  65 )(67 65)(70 60)

\move(65 62) \textref h:R v:C \htext{\tiny{\bf 1}}
\move(65 67) \textref h:R v:C \htext{\tiny{\bf 2}}
\move(45 62) \textref h:R v:C \htext{\tiny{\bf 3}}
\move(45 67) \textref h:R v:C \htext{\tiny{\bf 4}}

\move(80 57)\textref h:R v:C \htext{\footnotesize{2}}
\move(71 57)\textref h:R v:C \htext{\footnotesize{12}}
\move(61 57)\textref h:R v:C \htext{\footnotesize{\underline{4}2}}
\move(51 57)\textref h:R v:C \htext{\footnotesize{34}}

\move(30 45) \fcir f:0 r:0.4
\move(40 45) \fcir f:0 r:0.4
\move(50 45) \fcir f:0 r:0.4
\move(60 45) \fcir f:0 r:0.4
\move(70 45) \fcir f:0 r:0.4
\move(80 45) \fcir f:0 r:0.4

\move(30 45) \clvec (42  52 )(48 52)(60 45)
\move(40 45) \clvec (43  50 )(47 50)(50 45)
\move(50 45) \clvec (58  51 )(62 51)(70 45)
\move(60 45) \clvec (68  51 )(72 51)(80 45)

\move(70 51) \textref h:R v:C \htext{\tiny{\bf 1}}
\move(60 51) \textref h:R v:C \htext{\tiny{\bf 2}}
\move(45 47) \textref h:R v:C \htext{\tiny{\bf 3}}
\move(45 52) \textref h:R v:C \htext{\tiny{\bf 4}}

\move(80 42)\textref h:R v:C \htext{\footnotesize{1}}
\move(71 42)\textref h:R v:C \htext{\footnotesize{\underline{2}1}}
\move(61 42)\textref h:R v:C \htext{\footnotesize{\underline{4}2}}
\move(51 42)\textref h:R v:C \htext{\footnotesize{34}}

\move(110 90) \fcir f:0 r:0.4
\move(120 90) \fcir f:0 r:0.4
\move(130 90) \fcir f:0 r:0.4
\move(140 90) \fcir f:0 r:0.4
\move(150 90) \fcir f:0 r:0.4
\move(160 90) \fcir f:0 r:0.4

\move(110 90) \clvec (118  96 )(122 96)(130 90)
\move(120 90) \clvec (135  97 )(145 97)(160 90)
\move(130 90) \clvec (133  95 )(137 95)(140 90)
\move(140 90) \clvec (143  95 )(147 95)(150 90)

\move(145 92) \textref h:R v:C \htext{\tiny{\bf 1}}
\move(135 92) \textref h:R v:C \htext{\tiny{\bf 2}}
\move(140 96) \textref h:R v:C \htext{\tiny{\bf 3}}
\move(120 96) \textref h:R v:C \htext{\tiny{\bf 4}}

\move(160 87)\textref h:R v:C \htext{\footnotesize{3}}
\move(151 87)\textref h:R v:C \htext{\footnotesize{13}}
\move(141 87)\textref h:R v:C \htext{\footnotesize{23}}
\move(131 87)\textref h:R v:C \htext{\footnotesize{\underline{4}3}}

\move(110 75) \fcir f:0 r:0.4
\move(120 75) \fcir f:0 r:0.4
\move(130 75) \fcir f:0 r:0.4
\move(140 75) \fcir f:0 r:0.4
\move(150 75) \fcir f:0 r:0.4
\move(160 75) \fcir f:0 r:0.4

\move(110 75) \clvec (118  81)(122 81)(130 75)
\move(120 75) \clvec (132  82)(137 82)(150 75)
\move(130 75) \clvec (133  80)(137 80)(140 75)
\move(140 75) \clvec (148  81)(152 81)(160 75)

\move(150 81) \textref h:R v:C \htext{\tiny{\bf 1}}
\move(135 77) \textref h:R v:C \htext{\tiny{\bf 2}}
\move(135 81) \textref h:R v:C \htext{\tiny{\bf 3}}
\move(120 81) \textref h:R v:C \htext{\tiny{\bf 4}}

\move(160 72)\textref h:R v:C \htext{\footnotesize{1}}
\move(151 72)\textref h:R v:C \htext{\footnotesize{\underline{3}1}}
\move(141 72)\textref h:R v:C \htext{\footnotesize{23}}
\move(131 72)\textref h:R v:C \htext{\footnotesize{\underline{4}3}}

\move(110 60) \fcir f:0 r:0.4
\move(120 60) \fcir f:0 r:0.4
\move(130 60) \fcir f:0 r:0.4
\move(140 60) \fcir f:0 r:0.4
\move(150 60) \fcir f:0 r:0.4
\move(160 60) \fcir f:0 r:0.4

\move(110 60) \clvec (118  66 )(122 66)(130 60)
\move(120 60) \clvec (128  66 )(132 66)(140 60)
\move(130 60) \clvec (142  67 )(148 67)(160 60)
\move(140 60) \clvec (143  65 )(147 65)(150 60)

\move(145 62) \textref h:R v:C \htext{\tiny{\bf 1}}
\move(145 67) \textref h:R v:C \htext{\tiny{\bf 2}}
\move(130 66) \textref h:R v:C \htext{\tiny{\bf 3}}
\move(120 66) \textref h:R v:C \htext{\tiny{\bf 4}}

\move(160 57)\textref h:R v:C \htext{\footnotesize{2}}
\move(151 57)\textref h:R v:C \htext{\footnotesize{12}}
\move(141 57)\textref h:R v:C \htext{\footnotesize{\underline{3}2}}
\move(131 57)\textref h:R v:C \htext{\footnotesize{\underline{4}3}}

\move(110 45) \fcir f:0 r:0.4
\move(120 45) \fcir f:0 r:0.4
\move(130 45) \fcir f:0 r:0.4
\move(140 45) \fcir f:0 r:0.4
\move(150 45) \fcir f:0 r:0.4
\move(160 45) \fcir f:0 r:0.4

\move(110 45) \clvec (117  50 )(123 50)(130 45)
\move(120 45) \clvec (127  50 )(133 50)(140 45)
\move(130 45) \clvec (137  50 )(143 50)(150 45)
\move(140 45) \clvec (147  50 )(153 50)(160 45)

\move(150 50) \textref h:R v:C \htext{\tiny{\bf 1}}
\move(140 50) \textref h:R v:C \htext{\tiny{\bf 2}}
\move(130 50) \textref h:R v:C \htext{\tiny{\bf 3}}
\move(120 50) \textref h:R v:C \htext{\tiny{\bf 4}}

\move(160 42)\textref h:R v:C \htext{\footnotesize{1}}
\move(151 42)\textref h:R v:C \htext{\footnotesize{\underline{2}1}}
\move(141 42)\textref h:R v:C \htext{\footnotesize{\underline{3}2}}
\move(131 42)\textref h:R v:C \htext{\footnotesize{\underline{4}3}}
}
\caption{$\pmaj(P)$ for $P \in P_6(\{1,2\}, \{5,6\})$. }
\label{major-ex}
\end{figure}
\end{ex}

\begin{remark}
{\rm The joint distribution of $\Cr_2$ and $\pmaj$ is not in
general symmetric over $P_n(S, T)$. For example,  let  $n=8$,
$S=\{1,2,3,5\}$, and $T=\{4, 6, 7,8\}$. Then
\begin{multline*}
\sum_{P \in P_8(S,T)} x^{\Cr_2(P)} y^{\pmaj(P)} =
x^5y^5+x^4y^4+2x^3y^4+2x^4y^3 \\ +x^3y^3+3x^2y^2+2xy+1+2x^3y^2+x^2y^3+x^2y+xy^3.
\end{multline*} }
\end{remark}
\begin{remark} {\rm
 We explain the combinatorial meaning of the quantities
$\{\, h(i)=|T\ri|-|S\ri|: i \notin T\,\}$.
The paper \cite{CDDSY}   gives a characterization
of nonempty   $P_n(S,T)$'s.
Given a pair $(S,T)$ where $S, T \subseteq [n]$ and $|S|=|T|$,
associate to it a lattice path
$L(S, T)$  with steps $(1, 1)$, $(1, -1)$ and $(1, 0)$:
start from $(0,0)$, read the integers $i$ from $1$ to $n$ one by one, and
move  two steps for each $i$ as follows. \\
\mbox{} \hspace{2em} 1. If $i \in S \cap T$, move $(1, 0)$ twice.\\
\mbox{} \hspace{2em} 2. If $i \in S\setminus T$, move $(1,0)$ and
then $(1,1)$.\\
\mbox{} \hspace{2em} 3. If $i \in T\setminus S$, move $(1,-1)$ and
then $(1,0)$.\\
\mbox{} \hspace{2em} 4. If $i \notin S\cup T$, move $(1,-1)$ and then
$(1,1)$. \\
This defines a lattice path $L(S, T)$ from $(0,0)$ to $(2n, 0)$.
Conversely, the path uniquely determines $(S, T)$.
Then $P_n(S, T) $ is nonempty if and only if the lattice path $L(S,
T)$ is a Motzkin path, i.e., never goes below the $x$-axis.

For each element $i \in [n]\setminus T$, there is a unique upper step $(1,1)$
in the lattice path $L(S,T)$.
We say  an upper step is \emph{of height $y$} if it goes from $(x-1,y-1)$ to
$(x, y)$.
Then the multiset  $\{\, h(i)=|T\ri|-|S\ri|: i \notin T\,\}$
is exactly the same as the multiset
\{\,height of  $U$: $U$ is an upper step in  $L(S, T)\,\}$.
}
\end{remark}
%

\section{The generalized Foata bijection}

In this section we construct a bijection $\phi$ from $P_n(S,\,T)$
to itself such that $\pmaj(P)=\Cr_2(\phi(P))$ for any set partition $P$. This
provides a generalization for Foata's second fundamental
transformation \cite{Foata,Foata97} which is used to
prove the equidistribution of the permutation statistics
$\inv$ and $\maj$.

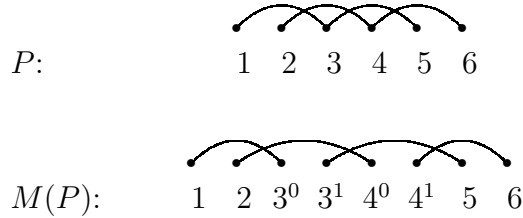
\begin{figure}[h,t]
\begin{center}
\begin{picture}(350,100)
\setlength{\unitlength}{6mm}

\multiput(8,4)(1,0){6}{\circle*{0.2}}

\qbezier(8,4)(9,5)(10,4) \qbezier(10,4)(11,5)(12,4)
\qbezier(9,4)(10,5)(11,4) \qbezier(11,4)(12,5)(13,4)

\put(3,3){$P$:}
\put(8,3){1} \put(9,3){2} \put(10,3){3} \put(11,3){4}
\put(12,3){5} \put(13,3){6}

\multiput(7,1)(1,0){8}{\circle*{0.2}}

\qbezier(7,1)(8,2)(9,1) \qbezier(10,1)(11.5,2)(13,1)
\qbezier(8,1)(9.5,2)(11,1) \qbezier(12,1)(13,2)(14,1)

\put(3,0){$M(P)$:}
\put(7,0){1} \put(8,0){2} \put(8.8,0){$3^0$} \put(9.8,0){$3^1$}
\put(10.8,0){$4^0$} \put(11.8,0){$4^1$} \put(13,0){5}
\put(14,0){6}

\end{picture}
\end{center}
\caption{A set partition $P$  and the corresponding matching $M(P)$.
}\label{fig-p-m}
\end{figure}

Given a partition $P$, for each endpoint $i\not\in S\cup T$, we may
replace $i$ with two neighboring endpoints, i.e., a right-hand
endpoint $i^0$ on the left and a left-hand endpoint $i^1$ on the
right, such that the  arc ending at $i$ is incident to $i^0$ and the
arc starting from $i$ is incident to $i^1$. After dealing with each
endpoint not in $S\cup T$ and removing all isolated points from $P$,
we obtain a matching $M(P)$. Clearly $\pmaj(P)=\pmaj(M(P))$ and
$\Cr_2(P)=\Cr_2(M(P))$. See Figure \ref{fig-p-m} for an example. Then
it is sufficient to describe the bijection $\phi$ on matchings, as
for any set partition $P$, one can obtain $\phi(P)$ from
$\phi(M(P))$ by  identifying $i^0$ and $i^1$ as one endpoint and
adding back all isolated vertices.

We introduce some notation.  Let $M$ be a matching of $[2m]$.
Suppose the arc incident to  $2m$ is $(i, 2m)$.
We say an arc $(a,b)$ of $M$ is {\em large} if $a < i$,
and {\em small} if $a > i$.
This terminology can also be described in terms of the edge labeling.
Recall that
if $(i_1, j_1), (i_2, j_2), \dots$, $(i_m,j_m)$ are the arcs of $M$
where $i_1 > i_2 > \cdots >i_m$, then the arc $(i_k, j_k)$ has
label $k$.
Suppose the arc $(i, 2m)$ has label $b$. Then an arc is large if
its label is larger than $b$, and small if its label is smaller
than  $b$.

Given a matching $M$ on $[2m]$ with the arc $(i, 2m)$, let $N = M
\setminus (i, 2m)$. Let $\mathcal{R}$ be the set of right-hand
endpoints of $N$ lying between $i$ and $2m$. We divide the set
$\mathcal{R}$ into three disjoint subsets. Define the {\em critical
large arc} $e_L$ to be the large arc with the biggest right-hand
endpoint in $\mathcal{R}$, and the {\em critical small arc} $e_S$ to
be the small arc with the smallest left-hand endpoint which crosses
$e_L$. Assume $e_S=(e_S^l, e_S^r)$. Then we set
\begin{align*}
\mR_0(N)&=\{\, j \in \mR:  e_S^r < j < 2m\,\}, \\
\mR_1(N)&=\{\,  j \in \mR: e_S^l < j \leq e_S^r\,\}, \\
\mR_2(N)&=\{\, j \in  \mR:  i < j < e_S^l\,\}.
\end{align*}
If there exists no critical large arc $e_L$, then $\mR_1(N)$ and
$\mR_2(N)$
are empty, and all endpoints belong to $\mR_0(N)$.
If $e_L=(e_L^l, e_L^r)$ exists
but no small arc crosses it, then
$\mR_2(N)$ is composed of all right-hand endpoints $j$ with $i < j \leq
e_L^r$, $\mR_1(N)$ is empty, and $\mR_0(N)$ is composed of the remaining
endpoints in $\mR$.

\begin{ex}
Let $M$ be the matching with arcs $(1, 14)$, $(2, 7)$, $(3, 16)$, $(4,6)$,
$(5,9)$, $(8, 12)$, $(10,15)$ and $(11,13)$.
Let $N=M \setminus (3,16)$. In $N$
the arcs $(1, 14)$ and $(2, 7)$ are large arcs.  $(4, 6), (5, 9), (8, 12),
(10,15)$ and $(11, 13)$ are small arcs.
 The critical large arc is $e_L=(1,14)$,
and  the  critical small arc $e_S=(10,15)$. The set $\mR$ consists
of vertices $\{6, 7, 9, 12, 13, 14, 15\}$, where $\mR_0=\emptyset$,
$\mR_1=\{12, 13, 14, 15\}$, and $\mR_2=\{6, 7, 9\}$.  See Figure
\ref{constru012}.

\begin{figure}[h,t]
\begin{center}
\begin{picture}(300,40)
\setlength{\unitlength}{3.5mm}

\multiput(0,0)(2,0){16}{\circle*{0.2}} \put(-0.3,-1){1}
\put(1.7,-1){2}\put(3.7,-1){3}  \put(5.7,-1){4} \put(7.7,-1){5}
\put(9.7,-1){6} \put(11.7,-1){7} \put(13.7,-1){8} \put(15.7,-1){9}
\put(17.5,-1){10} \put(19.5,-1){11} \put(21.5,-1){12}
\put(23.5,-1){13} \put(25.5,-1){14} \put(27.5,-1){15}
\put(29.5,-1){16}

\qbezier(2,0)(7, 3)(12,0) \put(6, 1.5){\tiny{7}}

\qbezier(4,0)(17, 6)(30,0)
\put(22, 3.1){\tiny{6}}

\qbezier(6,0)(8,1)(10,0)
\put(8,0.8){\tiny{5}}

\qbezier(8,0)(12,2)(16,0)
\put(12,1.3){\tiny{4}}

\qbezier(14,0)(18,2)(22,0)
\put(18,1.4){\tiny{3}}

\qbezier(20,0)(22,1.5)(24,0)
\put(22,0.8){\tiny{1}}

\put(23,1.6){\tiny{2}} \put(12, 3.3){\tiny{8}}


\linethickness{1pt}

\qbezier[40](18,0)(23,3)(28,0) \qbezier[60](0,0)(14,6)(26,0)

\end{picture}
\end{center}
\caption{The construction of $\mR_0$, $\mR_1$ and
$\mR_2$.}\label{constru012}
\end{figure}
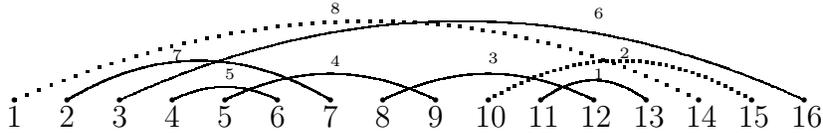
\end{ex}

In the following we describe the map $\phi$ on the set of matchings,
that is, on $P_n(S, T)$ where $n=2m$, $|S|=|T|=m$, and $S \cap
T=\emptyset$. The map $\phi$ preserves
the arc incident to $2m$;  i.e., if $(i, 2m)$ is an arc of
$M$, then it is also an arc of $\phi(M)$.
We extend the map $\phi$ to matchings whose vertices
are $a_1 < a_2 < \cdots < a_{2m}$ by identifying the vertex $a_i$
with $i$.

The map  $\phi$ is defined by induction on  $m$.
When $m=1$, let $\phi$ be the identity map. Given a
matching  $M$ of $2m$ for $m>1$,
let $M_1$ be the matching obtained from $M$ by removing the
arc $(i, 2m)$, and $M_1'=\phi(M_1)$. We construct
a matching $M_2$  from $M_1'$ by applying
a series of operations  on the arcs of $M_1'$ whose
right-hand endpoints are in the set $\mR$.

\medskip

\noindent {\bf 1. On $\mR_0(M_1')$.}

 We fix all the arcs whose right-hand endpoints are in $\mR_0(M_1')$.

\medskip
\noindent {\bf 2. Algorithm I on $\mR_1(M_1')$.}

Let $ptr_1$, $ptr_2$ be two pointers. We apply the following
algorithm on  $\mR_1(M_1')$ if $\mR_1(M_1')$ is nonempty.

{\sf
\begin{itemize}
\item[{\bf(A)}] Let $ptr_1$ point to $e_S^r$,  and  $ptr_2$ point to
the next vertex in $\mR_1(M_1')$ on the left of $e_S^r$.

\item[{\bf(B)}] If $ptr_2$ is null, then go to (D). Otherwise,
assume $ptr_1=j_1$ and $ptr_2=j_2$, where $j_1, j_2$ are right-hand
endpoints of the arcs $(i_1, j_1)$ and $(i_2, j_2)$.

\begin{itemize}

   \item[{\bf(B1)}] If $(i_2, j_2)$ is a large arc,
   then change the two arcs $\{(i_1, j_1), (i_2, j_2)\}$  to
   $\{(i_1,j_2), (i_2,j_1)\}$. Move $ptr_1$  to $j_2$.

   \item[\bf{(B2)}] If $(i_2, j_2)$ is  a small arc $(i_2,j_2)$,
        then there are three cases to consider.

   \begin{itemize}

   \item[{\bf(B2.1)}] If there exists no right-hand endpoint between
    $i_1$ and $i_2$, then  move $ptr_1$  to $j_2$.

   \item[{\bf(B2.2)}] If $i_1 < i_2$ and there are some right-hand endpoints
 between  $i_1$ and $i_2$, do nothing.

   \item[{\bf(B2.3)}] If $i_2 < i_1$ and there are some right-hand
   endpoints
   between  $i_1$ and $i_2$, then find the smallest $j_3$ such that
  $j_3 > j_1$ and   $(i_3,j_3)$ is a large arc.
    Change the three arcs
   $\{(i_1,j_1), (i_2, j_2), (i_3, j_3)\}$ to
  $\{(i_1,j_2), (i_2,j_3), (i_3,j_1)\}$, and move $ptr_1$  to $j_2$.
  See Figure \ref{b2.3} for an illustration, where $j$ is a right-hand endpoint.

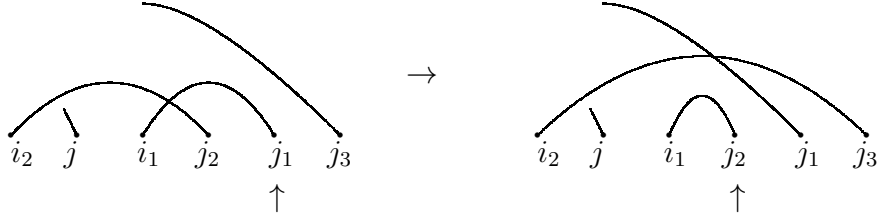
\begin{figure}[h,t]
\begin{center}
\begin{picture}(350,75)
\setlength{\unitlength}{3.5mm}

\multiput(2,2)(2.5,0){6}{\circle*{0.2}}
\multiput(22,2)(2.5,0){6}{\circle*{0.2}}

\qbezier(2,2)(5.75,6)(9.5,2) \qbezier(22,2)(28.25,8)(34.5,2)
\qbezier(7,2)(9.5,6)(12,2) \qbezier(27,2)(28.25,5)(29.5,2)
\qbezier(4,3)(4.25, 2.5)(4.5,2) \qbezier(7,7)(9.5,7)(14.5,2)
\qbezier(24,3)(24.25,2.5)(24.5,2) \qbezier(24.5,7)(27,7)(32,2)

\put(2,1.0){{$i_2$}} \put(9,1.0){{$j_2$}} \put(6.8,1.0){{$i_1$}}
\put(11.8,1.0){{$j_1$}} \put(14,1.0){{$j_3$}}

\put(22,1.0){{$i_2$}} \put(29,1.0){{$j_2$}} \put(26.8,1.0){{$i_1$}}
\put(31.8,1.0){{$j_1$}} \put(34,1.0){{$j_3$}}

\put(29.3,-0.7){{$\uparrow$}} \put(11.8,-0.7){{$\uparrow$}}

\put(4,1.0){$j$} \put(24,1.0){$j$}

\put(17,4){$\rightarrow$}

\end{picture}
\end{center}
\caption{Step (B2.3), where the upper arrow represents the pointer
$ptr_1$.}\label{b2.3}
\end{figure}

\end{itemize}
\end{itemize}

\item[{\bf(C)}] Move $ptr_2$ to the next vertex in $\mR_1(M_1')$
on the left of $j_2$. Go to (B).

\item[{\bf(D)}] Stop.
\end{itemize}
}

\noindent {\bf 3. Algorithm II on $\mR_2(M_1')$.}

We continue by processing vertices in  $\mR_2(M_1')$. Let the
rightmost vertex of $\mR_2(M_1')$ be~$t_0$, which is the right-hand
endpoint of the arc $(i_0, t_0)$.
{\sf
\begin{itemize}
\item ({\bf Pre-process}):  If $(i_0, t_0)$ is a small arc,
then find the large arc $(i_1,t_1)$ on the right of
$(i_0,t_0)$ with smallest $t_1$.
By the definition of $\mR_1$ and $\mR_2$, such a large
arc exists. We change the two arcs $\{(i_0, t_0), (i_1,t_1)\}$
to $\{(i_0,t_1), (i_1,t_0)\}$.
Now the vertex $t_0$ is connected to a large arc.
\end{itemize}
}
\begin{remark} \label{borrowing}
 We call this step {\em borrowing}. Note that now
the arc $(i_0, t_1)$ is a small arc which crosses $(i_1, t_0)$, but
does not cross any other large arc in $\mR_1$, and $t_1$ is the smallest
vertex in $\mR_1$ with this property.
\end{remark}

Let $ptr_1, ptr_2$ be two pointers.
We apply the following operations on  $\mR_2(M_1')$.

{\sf
\begin{itemize}
\item[{\bf(A$'$)}] Let $ptr_1$ point to $t_0$, and let $ptr_2$ point to the next
vertex in $\mR$ on the left of $t_0$.

\item[{\bf(B$'$)}] If $ptr_2$ is null, then go to (D$'$). Otherwise,
assume $ptr_1=j_1$ and $ptr_2=j_2$, where $j_1, j_2$ are right-hand
endpoints of the arcs $(i_1, j_1)$ and $(i_2, j_2)$.

\begin{itemize}

   \item[{\bf(B$'1$)}] If $(i_2, j_2)$ is a large arc, then do nothing.

   \item[{\bf(B$'2$)}] If $(i_2, j_2)$ is a small arc, then change the two
 arcs  to  $(i_1,j_2)$ and $(i_2,j_1)$.

\end{itemize}

\item[{\bf(C$'$)}]  Move $ptr_1$ to $j_2$, and move $ptr_2$  to
the next vertex in $\mR$ on the left of $j_2$. Go to (B$'$).

\item[{\bf(D$'$)}] Stop.
\end{itemize}
}

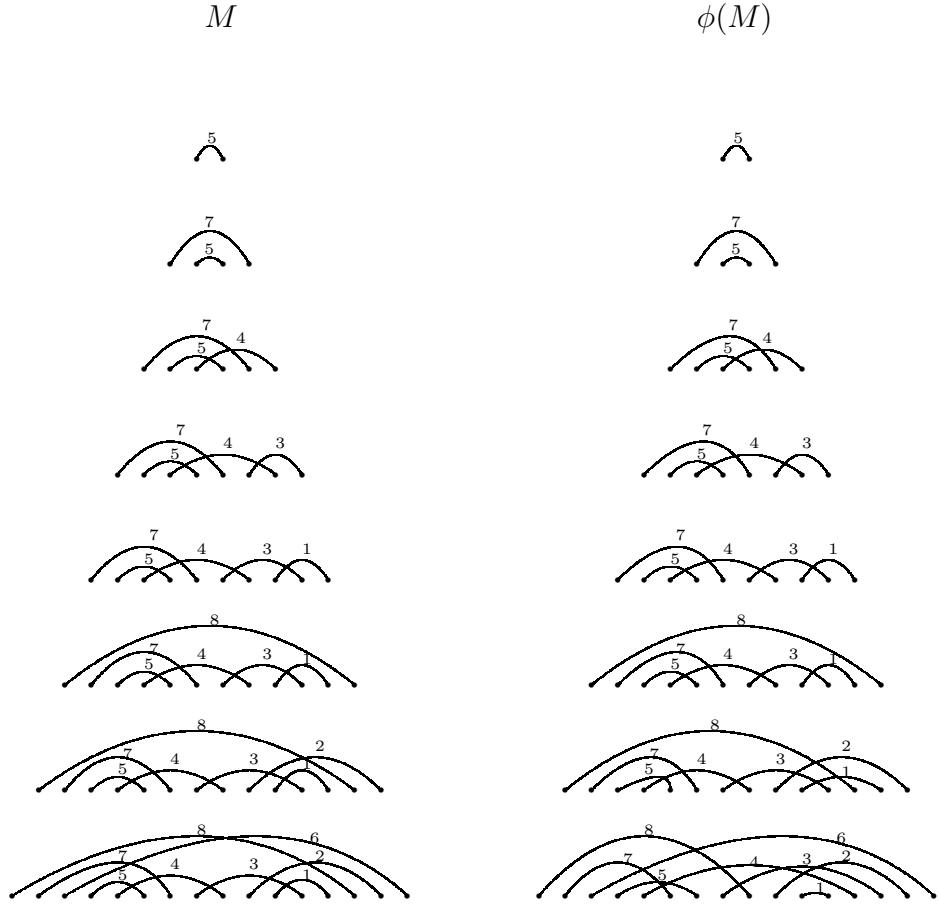
\begin{figure}[h,t]
\begin{center}
\begin{picture}(350,330)
\setlength{\unitlength}{3.5mm}

\put(7.3,33){$M$} \put(26,33){$\phi(M)$}

\multiput(7,28)(1,0){2}{\circle*{0.2}}
\multiput(27,28)(1,0){2}{\circle*{0.2}}

\qbezier(7,28)(7.5,29)(8,28) \qbezier(27,28)(27.5,29)(28,28)
\put(7.4,28.6){\tiny{5}} \put(27.4,28.6){\tiny{5}}


\multiput(6,24)(1,0){4}{\circle*{0.2}}
\multiput(26,24)(1,0){4}{\circle*{0.2}}

\qbezier(7,24)(7.5,24.5)(8,24) \qbezier(27,24)(27.5,24.5)(28,24)
\qbezier(6,24)(7.5,26.5)(9,24) \qbezier(26,24)(27.5,26.5)(29,24)
\put(7.3,24.4){\tiny{5}} \put(27.3,24.4){\tiny{5}}
\put(7.3,25.4){\tiny{7}} \put(27.3,25.4){\tiny{7}}

\multiput(5,20)(1,0){6}{\circle*{0.2}}
\multiput(25,20)(1,0){6}{\circle*{0.2}}

\qbezier(5,20)(7,22.5)(9,20) \qbezier(25,20)(27.5,22.5)(29,20)
\put(7.2,21.5){\tiny{7}} \put(27.2,21.5){\tiny{7}}

\qbezier(6,20)(7,21)(8,20) \qbezier(26,20)(27,21)(28,20)
\put(7,20.6){\tiny{5}} \put(27,20.6){\tiny{5}}

\qbezier(7,20)(8.5,21.5)(10,20) \qbezier(27,20)(28.5,21.5)(30,20)
\put(8.5,21){\tiny{4}} \put(28.5,21){\tiny{4}}

\multiput(4,16)(1,0){8}{\circle*{0.2}}
\multiput(24,16)(1,0){8}{\circle*{0.2}}

\qbezier(4,16)(6,18.5)(8,16) \qbezier(24,16)(26.5,18.5)(28,16)
\put(6.2,17.5){\tiny{7}} \put(26.2,17.5){\tiny{7}}

\qbezier(5,16)(6,17)(7,16) \qbezier(25,16)(26,17)(27,16)
\put(6,16.6){\tiny{5}} \put(26,16.6){\tiny{5}}

\qbezier(6,16)(8,17.5)(10,16) \qbezier(26,16)(28,17.5)(30,16)
\put(8,17){\tiny{4}} \put(28,17){\tiny{4}}

\qbezier(9,16)(10,17.5)(11,16) \qbezier(29,16)(30,17.5)(31,16)
\put(10,17){\tiny{3}} \put(30,17){\tiny{3}}

\multiput(3,12)(1,0){10}{\circle*{0.2}}
\multiput(23,12)(1,0){10}{\circle*{0.2}}

\qbezier(3,12)(5,14.5)(7,12) \qbezier(23,12)(25.5,14.5)(27,12)
\put(5.2,13.5){\tiny{7}} \put(25.2,13.5){\tiny{7}}

\qbezier(4,12)(5,13)(6,12) \qbezier(24,12)(25,13)(26,12)
\put(5,12.6){\tiny{5}} \put(25,12.6){\tiny{5}}

\qbezier(5,12)(7,13.5)(9,12) \qbezier(25,12)(27,13.5)(29,12)
\put(7,13){\tiny{4}} \put(27,13){\tiny{4}}

\qbezier(8,12)(9.5,13.5)(11,12) \qbezier(28,12)(29.5,13.5)(31,12)
\put(9.5,13){\tiny{3}} \put(29.5,13){\tiny{3}}

\qbezier(10,12)(11,13.5)(12,12) \qbezier(30,12)(31,13.5)(32,12)
\put(11,13){\tiny{1}} \put(31,13){\tiny{1}}

\multiput(2,8)(1,0){12}{\circle*{0.2}}
\multiput(22,8)(1,0){12}{\circle*{0.2}}

\qbezier(2,8)(7.5,12.5)(13,8) \qbezier(22,8)(27.5,12.5)(33,8)
\put(7.5, 10.3){\tiny{8}} \put(27.5,10.3){\tiny{8}}

\qbezier(3,8)(5,10.5)(7,8) \qbezier(23,8)(25.5,10.5)(27,8)
\put(5.2,9.2){\tiny{7}} \put(25.2,9.2){\tiny{7}}

\qbezier(4,8)(5,9)(6,8) \qbezier(24,8)(25,9)(26,8)
\put(5,8.6){\tiny{5}} \put(25,8.6){\tiny{5}}

\qbezier(5,8)(7,9.5)(9,8) \qbezier(25,8)(27,9.5)(29,8)
\put(7,9){\tiny{4}} \put(27,9){\tiny{4}}

\qbezier(8,8)(9.5,9.5)(11,8) \qbezier(28,8)(29.5,9.5)(31,8)
\put(9.5,9){\tiny{3}} \put(29.5,9){\tiny{3}}

\qbezier(10,8)(11,9.5)(12,8) \qbezier(30,8)(31,9.5)(32,8)
\put(11,8.8){\tiny{1}} \put(31,8.8){\tiny{1}}

\multiput(1,4)(1,0){14}{\circle*{0.2}}
\multiput(21,4)(1,0){14}{\circle*{0.2}}

\qbezier(1,4)(7,8.5)(13,4) \qbezier(21,4)(26.5,8.5)(32,4) \put(7,
6.3){\tiny{8}} \put(26.5,6.3){\tiny{8}}

\qbezier(2,4)(4,6.5)(6,4) \qbezier(22,4)(24.5,6.5)(26,4)
\put(4.2,5.2){\tiny{7}} \put(24.2,5.2){\tiny{7}}

\qbezier(3,4)(4,5)(5,4) \qbezier(23,4)(25,5)(25,4)
\put(4,4.6){\tiny{5}} \put(24,4.6){\tiny{5}}

\qbezier(4,4)(6,5.5)(8,4) \qbezier(24,4)(26,5.5)(28,4)
\put(6,5){\tiny{4}} \put(26,5){\tiny{4}}

\qbezier(7,4)(9,5.5)(11,4) \qbezier(27,4)(29,5.5)(31,4)
\put(9,5){\tiny{3}} \put(29,5){\tiny{3}}

\qbezier(10,4)(11,5.5)(12,4) \qbezier(30,4)(31.5,5)(33,4)
\put(11,4.8){\tiny{1}} \put(31.5,4.5){\tiny{1}}

\qbezier(9,4)(11.5,6.5)(14,4) \qbezier(29,4)(31.5,6.5)(34,4)
\put(11.5,5.5){\tiny{2}} \put(31.5,5.5){\tiny{2}}

\multiput(0,0)(1,0){16}{\circle*{0.2}}
\multiput(20,0)(1,0){16}{\circle*{0.2}}

\qbezier(0,0)(7.5,4.5)(13,0) \qbezier(20,0)(24,4.5)(28,0)
\put(7,2.3){\tiny{8}} \put(24,2.3){\tiny{8}}

\qbezier(1,0)(4.5,2.5)(6,0) \qbezier(21,0)(23,2.5)(25,0)
\put(4,1.3){\tiny{7}} \put(23.2,1.2){\tiny{7}}

\qbezier(2,0)(10,4.5)(15,0) \qbezier(22,0)(30,4.5)(35,0)
\put(11.3,2){\tiny{6}} \put(31.3,2){\tiny{6}}

\qbezier(3,0)(4,1)(5,0) \qbezier(23,0)(24.5,1)(26,0)
\put(4,0.6){\tiny{5}} \put(24.5,0.6){\tiny{5}}

\qbezier(4,0)(6,1.5)(8,0) \qbezier(24,0)(28,2.3)(32,0)
\put(6,1){\tiny{4}} \put(28,1.1){\tiny{4}}

\qbezier(7,0)(9,1.5)(11,0) \qbezier(27,0)(30,2.2)(33,0)
\put(9,1){\tiny{3}} \put(30,1.2){\tiny{3}}

\qbezier(10,0)(11,1.2)(12,0) \qbezier(30,0)(30.5,0.2)(31,0)
\put(11,0.6){\tiny{1}} \put(30.5,0.1){\tiny{1}}

\qbezier(9,0)(11.5,2.5)(14,0) \qbezier(29,0)(31.5,2.5)(34,0)
\put(11.5,1.3){\tiny{2}} \put(31.5,1.3){\tiny{2}}

\end{picture}
\end{center}
\caption{An example of the bijection $\phi$ }\label{bijection}
\end{figure}

Let $M_2$ be the matching obtained by applying {Algorithms I and II}
on $M_1'$. Then $\phi(M)=M'=M_2 \cup \{(i, 2m)\}$. We give an
example to illustrate the bijection in Figure \ref{bijection}.

%

To see that $\phi(M)$ is well-defined, we only need to check that {
Step (B2.3)} is valid.  Note that the algorithms do not change the
relative positions of the large arcs. In {Algorithm I}, if
$ptr_1=j_1$, then (1) $j_1$ must be the right-hand endpoint of a
small arc $(i_1, j_1)$, (2)  $i_1$ is smaller than any vertex in
$\mR_1(M_1')$, and (3) any right-hand endpoint between $ptr_1$ and
$ptr_2$ must be connected to a small arc.  Thus if $ptr_1=j_1$,
$ptr_2=j_2$ with arcs $(i_1, j_1)$ and $(i_2, j_2)$, where $i_2 <
i_1$ and there are some right-hand endpoints between $i_1$ and $i_2$,
then the critical large arc must have been moved to the right of the
pointers. Hence the large arc described in {Step (B2.3)} exists.
This shows that the map $\phi$ is a well-defined map. It is also
clear that $\phi$ preserves the set of left-hand endpoints, as well
as the set of  right-hand endpoints.

From the above construction, we notice  the following properties of
{Algorithms I and  II} and the matching $M_2$.

\medskip
\noindent {\bf Properties}:
\begin{enumerate}
\item For each vertex $j \in \mR_2(M_1')$, if $j$ is the right-hand
  endpoint of a small arc $(i,j)$  in $M_2$, then there is a large arc
  in $M_2$   that crosses $(i,j)$.

\item If $\mR_2(M_1')$ is nonempty, then the smallest vertex in $\mR_2(M_1')$
must connect to a large arc in $M_2$.

\item In {Algorithm I}, the pointer  $ptr_1$ always points to the right-hand
  endpoint of a small arc. In {Algorithm II}, $ptr_1$ always points
  to the right-hand endpoint of a large arc.

\item In {Algorithm I},  if the pointer $ptr_1$ is
connected to the small arc $(i_1,j_1)$, then $i_1$ is smaller than
any vertex in $\mR_1(M_1')$, and larger than any vertex in
$\mR_2(M_1')$.

\item Assume that when {Algorithm I} stops, the pointer $ptr_1$ is at
vertex $j_1$, which is the right-hand endpoint of a small arc $(i_1,
j_1)$. Then $j_1$ must be the smallest vertex in $\mR_1(M_1')$.
Algorithm II does not change the arc $(i_1, j_1)$, and there is no
arc $(i_2, j_2)$ in $M_2$ such that $i_2 < i_1 < j_2 < j_1$.
\end{enumerate}

\begin{thm}
 The map $\phi$, when restricted to $P_n(S, T)$,  is a bijection.
\end{thm}
\proof
It is sufficient to describe  how to invert {Algorithms I and II} in
the definition of  $\phi$, on the set of matchings. Given a matching
of $[2m]$ with the last arc $(i, 2m)$, let $N=M \setminus \{(i,
2m)\}$ and let $\mR$ be the set of right-hand endpoints of $N$ lying
between $i$ and $2m$. First we need to determine the set $\mR_2$. If
there exists no large arc whose right-hand endpoint lies between $i$
and $2m$, then $\mR_2=\mR_1=\emptyset$, and $\mR_0=\mR$. Suppose
there are some large arcs in $\mR$. Find the small arc $f=(f^l,
f^r)$ with the smallest right-hand endpoint
 $f^r$ in $\mR$ such that there is no arc $(i_2, j_2)$
with $i_2 < f^l$ which crosses $f$.
 By Properties 1 and 5 above, $\mR_2=\{ j \in \mR:  j < f^l\}$.
If no such small arc exists, then  $\mR_2=\mR$.

\medskip

\noindent {\bf Algorithm III: Inverse of Algorithm II.}\\
To invert {Algorithm II} in $\mR_2$, we apply the same steps as in
{Algorithm II} except that {\sf
\begin{itemize}
\item[(1)] Initially, let $ptr_1$ point to the smallest vertex in
$\mR_2$, and $ptr_2$ point to the next vertex in $\mR_2$ to
the right of $ptr_1$.
\item[(2)] We apply the {Steps  (B$'$)}--{(D$'$)} on the vertices in $\mR_2$
from left to right, i.e., in {Step (C$'$)} we need to move $ptr_2$
to the next vertex in $\mR$ on the right of $j_2$.
\end{itemize}
}

Note that the pointer $ptr_1$ always points to  the right-hand
endpoint of a large arc. When  $ptr_1$ reaches the rightmost vertex
$j_1$ of $\mR_2$, with a current large arc $(i_1,j_1)$, we need to
determine whether this arc $(i_1,j_1)$ is borrowed from $\mR_1$. By
Remark \ref{borrowing}, this can be done as follows:

{\sf
\begin{itemize}
\item {\bf(Invert borrowing)}
Let  $(i_2,j_2)$ be the small arc with minimal $j_2$ such
that  $i_1 < i_2 < j_1 < j_2.$
\begin{itemize}

\item[(a)] If no such small arc exists, then there is no
borrowing.

\item[(b)] If there exists a right-hand endpoint of a large arc
between  $j_1$ and $j_2$, then there is no borrowing.

\item[(c)] If there exists no right-hand endpoint of any large arc
between $j_1$ and $j_2$, then the arc  $(i_1,j_1)$ is borrowed from
$\mR_1$.
To invert, change the two arcs $\{(i_1,j_1), (i_2, j_2)\}$
to  $\{(i_1,j_2), (i_2,j_1)\}$.
\end{itemize}
\end{itemize}
} At this stage,  if there is a small arc $(i', j')$ with $i' < j_1
< j'$, then the arc $(i', j')$ must cross some large arc whose
right-hand endpoint  is  in $\mR\setminus \mR_2$.

We continue by inverting {Algorithm I} in $\mR\setminus \mR_2$, when
$\mR_2 \neq \mR$. In the following, $j$ always represents a right-hand
endpoint in $\mR\setminus \mR_2$.

\medskip

\goodbreak
\noindent {\bf Algorithm IV:  Inverse of Algorithm I.}

{\sf
\begin{itemize}
\item[{\bf(IA)}] Let $ptr_1$ point to $f^r$,  the smallest vertex
in $\mR\setminus \mR_2$. Let  $ptr_2$ point to the next right-hand
endpoint on the  right of $f^r$.

\item[{\bf(IB)}] If $ptr_2$ is null, then go to (ID).
Otherwise, assume $ptr_1=j_1$, $ptr_2=j_2$ where $j_1, j_2$ are the
right-hand endpoints of arcs $(i_1, j_1)$ and $(i_2, j_2)$.

\begin{itemize}

   \item[{\bf(IB1)}] If $(i_2,j_2)$ is a large arc,
     we need to consider two cases.

   \begin{itemize}

   \item[{\bf(IB1.1)}] There exists a small arc $(i_3,j_3)$ such that (i) no
   large arc lies between $j_2$ and $j_3$, (ii) $i_3 < i_1$ and $j_2<j_3$, and
   (iii) there are some right-hand
   endpoints between $i_1$ and $i_3$. In this case, we choose such a small
   arc $(i_3,j_3)$ with minimal $j_3$, and then change the
    three arcs $\{(i_1, j_1), (i_2, j_2), (i_3, j_3)\}$
to $\{(i_1,j_2), (i_2,j_3), (i_3,j_1)\}$.
   Move $ptr_1$ to $j_2$.

   \item[{\bf(IB1.2)}] There exists no  small arc $(i_3, j_3)$ satisfying
     the conditions above.
   In this case, we change the two arcs $\{(i_1, j_1), (i_2, j_2)\}$ to
  $\{(i_1,j_2), (i_2,j_1)\}$. Move $ptr_1$ to $j_2$.

   \end{itemize}

   \item[{\bf(IB2)}] If $(i_2, j_2)$ is a small arc, there are also
   two cases.

   \begin{itemize}

   \item[{\bf(IB2.1)}] If  there exists no right-hand endpoint between
     $i_1$ and $i_2$, then move $ptr_1$ to $j_2$.

   \item[{\bf(IB2.2)}] If  there exist some right-hand
    endpoints between $i_1$ and $i_2$, then do nothing.

   \end{itemize}

\end{itemize}

\item[{\bf(IC)}] Let $ptr_2$ point to the next right-hand endpoint to the right of
$j_2$. Go to (IB).

\item[{\bf(ID)}] Stop.
\end{itemize}
}

In running the above algorithm, if $ptr_1=j_1$, $ptr_2=j_2$
with the arcs $(i_1, j_1)$ and $(i_2, j_2)$, then the following
properties hold:
\begin{enumerate}
\item $(i_1, j_1)$ is a small arc. Any right-hand vertex $j$ of $\mR$ with $j <
i_1$ is in $\mR_2$, and any right-hand vertex $j$ with $j > i_1$ is in
$\mR \setminus \mR_2$.
\item If $(i', j')$ is an arc with $j_1 < j' <j_2$, then there exists
a right-hand vertex $j$ such that $i_1 < j < i'$.
\end{enumerate}
The above steps enable us to get a matching $N'$ such that when
applying {Algorithms I and II} to $N'$, one gets $N$. To see this,
first assume $(i_2, j_2)$ is a large arc.
\begin{itemize}
\item [(1)] If there is an arc $(i_3, j_3)
$ as described in {Step (IB1.1)}, then the current configuration can
be obtained by applying {Algorithm I}, {Step (B2.3)} to arcs
$\{(i_1,j_2),(i_2,j_3),(i_3,j_1)\}$ with $ptr_1=j_2$ and
$ptr_2=j_3$. {Step (IB1.1)} reverses this operation.
\item [(2)]
If there is no arc $(i_3, j_3)$ as described in Step (IB1.1), then
the current configuration can be obtained by applying
 {Algorithm I}, {Step (B1)} to $\{(i_1, j_2), (i_2, j_1)\}$. Step
 (IB1.2)
reverses it.
\end{itemize}

In the case that $(i_2, j_2)$ is a small arc,
\begin{itemize}
\item [(1)] If there exists no right-hand
endpoint between $i_1$ and $i_2$, then the current configuration can
be obtained by applying {Algorithm I}, {Step (B2.1)}. {Step (IB2.1)}
reverses this operation.
\item [(2)]
If there exist some right-hand endpoints between $i_1$ and $i_2$,
then $i_1 < i_2$. Otherwise, assume $i_2 < i_1$. Then any right-hand
endpoint $j$ between $i_1$ and $i_2$ must be in $\mR_2$. But at the
stage when {Algorithm III} and the step of {Invert borrowing}
stop, there are some large arcs $(i', j')$ crossing $(i_2, j_2)$,
and $j' \in \mR\setminus \mR_2$. By our construction, the arc $(i_2,
j_2)$ should have been destroyed by an application of {Step
(IB1.1)}. Contradiction! Hence $i_1 < i_2$, and the current
configuration can be obtained by applying {Algorithm I}, {Step
(B2.2)}. {Step (IB2.2)} reverses it.
\end{itemize}
Since there is no large arc in $\mR_0$, applying {Algorithm IV} in
$\mR\setminus \mR_2=\mR_1 \cup \mR_0$ will not change the arcs in
$\mR_0$.

Let $N'$ be the matching obtained from $N$ by applying {Algorithms
III and IV}. The above argument shows that

$$N\,\raisebox{1.6ex}{$\underrightarrow{\mbox{Algorithms I and II}}$}\,N'.$$

 It follows that the map $\phi$ is surjective. As the set
$P_n(S, T)$ is finite, $\phi$ must be a bijection. \qed

We say an arc  $(i_1,j_1)$ of a matching $M$
is {\em maximal} if there
is no arc $(i_2,j_2)$ in $M$ such that $i_2 < i_1 < j_1 < j_2$.
For a maximal arc $e=(i_1,j_1)$, let
$$
t(e, M)=\min\{\,i' : \mbox{There is an arc $(i',j')$ such that
$i_1 < i' < j_1 < j'$}\,\}.
$$
If there is no such arc, let $t(e, M)=j_1+1$.

\begin{lemma} \label{regionII}
\begin{enumerate}
\item An arc $(i_1, j_1)$ is  maximal in $M$ if and only if $(i_1, k)$
  is maximal  in $M'=\phi(M)$, for some $k$.
\item Let $e=(i_1, j_1)$ be a maximal arc in $M$, and $e'=(i_1,
  k)$ be the corresponding maximal edge in $M'$. Then
  $t(e,M)=t(e',M')$.
\end{enumerate}
\end{lemma}
\proof Let $M$ be a matching of $[2m]$. We proceed by induction on
$m$. The case $m=1$ is trivial. Assume the claim holds for all
matchings on $2m-2$ linearly ordered vertices. Let $M$ be a matching
of $[2m]$ with the arc $(i, 2m)$. Let $M_1=M \setminus \{(i, 2m)\}$,
let $M_1'=\phi(M_1)$, let $M_2$ be the matching obtained from $M_1'$
by applying {Algorithms I and II}, and  let $M'=\phi(M)=M_2 \cup
\{(i, 2m)\}$. Clearly $e=(i, 2m)$ is a maximal arc in both $M$ and
$M'$ with $t(e,M)=t(e, M')=2m+1$. Otherwise, an arc $e=(i_1, j_1)$
with $i_1 \neq i$ is a maximal arc  of $M$ (respectively $M'$) if
and only if $i_1 < i$ and $(i_1,j_1)$ is a maximal arc of $M_1$
(respectively $M_2$). Denote by $e$, $f$, and $e'$ the arcs whose
left-hand endpoint is $i_1$ in $M_1$, $M_1'$ and $M_2$,
respectively. Then  by the inductive hypothesis,  $e$ is  maximal in
$M_1$ if and only if $f$ is maximal in $M_1'$, in which case  $t(e,
M_1)=t(f, M_1')$. Since the algorithms do not change the relative
positions of large arcs, this happens if and only if $e'$ is maximal
in $M_2$.

Assume $e=(i_1,j_1)$ is maximal in $M_1$ and  $x=t(e, M_1)=t(f,M_1')$.
{\setlength{\leftmargini}{3\parindent}
 \begin{itemize}
  \item[(Case 1)]
If $x < i$, then $t(e, M)=x$. The equation $t(f, M_1')=x$ implies that in
$M_1'$, among all arcs $(i,j)$ which cross $f$ with $i > i_1$, the
one with the smallest left-hand endpoint is a large arc. This
property is preserved by the algorithms, so $t(e', M_2)=x$, and
hence $t(e', M')=x$.
  \item[(Case 2)] If $x >i$, then $t(e, M)=i$, as the last arc $(i, 2m)$
crosses every large arc. Since $t(f, M_1')=x$, there is no large
arc in $M_1'$ that crosses $f$. Hence there is no large arc in $M_2$
that crosses $e'$. So again we have $t(e', M')=i$.
\end{itemize}
}
In both cases, $t(e, M)=t(e',M')$.
\qed

\begin{thm}
 We have $\pmaj(M)=\Cr_2(\phi(M))$ for all $M \in P_n(S,T)$.
\end{thm}
\proof
Again it is enough to prove the theorem for matchings of $[2m]$.
We proceed by induction on $m$. The theorem is clearly true for $m=1$.
Assume it is true for all matchings on $2m-2$ linearly ordered vertices.
Given a matching $M$ of $[2m]$ with the arc $(i, 2m)$,
let $M_1$, $M_1'$, $M_2$, and $M'$ be defined as in Lemma \ref{regionII}.
By the inductive hypothesis, $\pmaj(M_1)=\Cr_2(M_1')$. Hence
\begin{align}
\pmaj(M) & =\pmaj(M_1)+\# \{\,j: j \in \mR_2(M_1)\,\} \nonumber \\
        & = \Cr_2(M_1')+\#\{\,j: j \in \mR_2(M_1)\,\}.     \label{pmaj}
\end{align}
By our construction {Algorithm I} decreases the crossing number of
$M_1'$  by 1 for each right-hand endpoint of a large arc in
$\mR_1(M_1')$, and {Algorithm II} increases the crossing number by~1
for each right-hand endpoint of a small arc in $\mR_2(M_1')$. Hence
\begin{eqnarray*}
\Cr_2(M_2)= \Cr_2(M_1')&-&\#\{\,j \in \mR_1(M_1'): \mbox{$j$ belongs to
a large arc} \,\}  \\
    & + &  \#\{\,j \in\mR_2(M_1'):   \mbox{$j$ belongs to a  small arc} \,\}.
\end{eqnarray*}
Therefore
\begin{eqnarray}
\Cr_2(M')=\Cr_2(M_2)&+& \#\{\,j: j \in \mR: \mbox{$j$ belongs to a large
  arc}\,\}
  \nonumber \\
   =\Cr_2(M_1')&-&\#\{\,j \in \mR_1(M_1'): \mbox{$j$ belongs to
a large arc} \,\} \nonumber \\
    &+ &  \#\{\,j \in\mR_2(M_1'):   \mbox{$j$ belongs to a  small arc}
   \,\} \nonumber\\
    &+ & \# \{\,j: j \in \mR: \mbox{$j$ belongs to a large arc}\,\} \nonumber\\
  = \Cr_2(M_1') & + & \# \{\,j: j \in \mR_2(M_1')\,\} \label{cr2}
 \end{eqnarray}
Comparing Eqs. \eqref{pmaj} and \eqref{cr2},  we only need to show that
\begin{eqnarray} \label{mR_2}
\mR_2(M_1)=\mR_2(M_1').
\end{eqnarray}
But the critical large arc $e_L$  of $M_1$, if it exists,  must be a maximal
arc of $M_1$, and the left-hand endpoint of the critical small arc $e_S$, if it
exists, must be
$t(e_L, M_1)$. Hence the identity  \eqref{mR_2} follows from
Lemma \ref{regionII}.
\qed

Finally we explain how our construction extends
Foata's second fundamental transformation
on permutations, which can be described as follows. Let
$w=w_1w_2\cdots w_n$ be a word on $\mathbb{N}$ and let $a \notin
\{w_1, \dots, w_n\}$. If $w_n < a$, the $a$-factorization of $w$ is
$w=v_1b_1\cdots v_pb_p$, where each $b_i$ is a letter less than $a$,
and each $v_i$ is a word (possibly empty), all of whose letters are
greater than $a$. Similarly, if $w_n > a$, the $a$-factorization of
$w$ is $w=v_1b_1\cdots v_pb_p$, where each $b_i$ is a letter greater
than  $a$, and each  $v_i$ is a word (possibly empty), all of whose
letters are less than $a$. In each case we define
$$
\gamma_a(w)=b_1v_1\cdots b_pv_p.
$$
With the above notation, let $a=w_n$ and let $w'=w_1\cdots w_{n-1}$.
The second fundamental transformation $\Phi$ is defined recursively by
$\Phi(w)=w$ if $w$ has length $1$, and
$$
\Phi(w)=\gamma_a(\Phi(w'))a,
$$
if $w$ has length $n >1$. The map $\Phi$ has the property that
$\inv (\Phi(w)) =\maj(w)$.

For a permutation $\pi$ of length $m$,
our bijection $\phi$, when applied to the matching
$M_\pi=\{\,(m+1-\pi(i), i+m): 1 \leq i \leq m\,\}$, is essentially the same as
Foata's transformation $\Phi(\pi)$. Note that the last arc of $M_\pi$
corresponds to the last entry $\pi(m)$ of $\pi$, and the set $\mR$
consists of all right-hand endpoints except $2m$.   Then
\begin{enumerate}
\item
If $\pi(m-1) < \pi(m)$, then $\mR_2=\emptyset$. The map $\gamma_a$
in Foata's transformation is equivalent to {Algorithm I}, where
cases (B2.2) and (B2.3) will not happen.

\item
If  $\pi(m-1) > \pi(m)$, then   $\mR_0=\mR_1=\emptyset$. The map
$\gamma_a$ in Foata's transformation  is equivalent to {Algorithm
II}.
\end{enumerate}


\end{document}